\documentclass[12pt]{article}

\usepackage{epsfig}
\usepackage{graphicx}
\usepackage{amsmath,amsthm}
\usepackage{mathbbol}
\usepackage{amssymb}
\usepackage{latexsym}
\usepackage[round]{natbib}
\usepackage{float}
\usepackage{color}
\usepackage{soul}

\newtheorem{theorem}{Theorem}
\newtheorem{lemma}{Lemma}

\newtheorem{proposition}{Proposition}

\newcommand{\giso}{{{\tilde g}}}
\newcommand{\ciso}{{{\tilde c}}}
\newcommand{\af}{\alpha}

\newcommand{\dd}{{\, \mathrm{d}}}
\newcommand{\bX}{\mathbf{X}}
\newcommand{\bY}{\mathbf{Y}}
\newcommand{\by}{\mathbf{y}}

\newcommand{\bfx}{\mathbf{x}}
\newcommand{\xn}{{\{x_1,\ldots,x_n\}}}
\newcommand{\xnl}{{x_1,\ldots,x_n}}
\newcommand{\bXred}{\mathbf{X}^!_\xnl}

\newcommand{\R}{\mathbb R}
\newcommand{\E}{\mathrm E}
\newcommand{\PP}{\mathrm P}

\newcommand{\T}{{\mbox{\scriptsize \sf T}}}
\usepackage{authblk}

\begin{document}
 
\title{
Palm distributions  \\ for log Gaussian Cox processes} 
\author[1]{Jean-Fran{\c{c}}ois Coeurjolly}

 \affil[1]{Laboratory Jean Kuntzmann, Grenoble Alpes University, France, email: \texttt{Jean-Francois.Coeurjolly@upmf-grenoble.fr}.}
\author[2]{Jesper M{\o}ller}
\author[2]{Rasmus Waagepetersen}  
\affil[2]{Department of Mathematical Sciences,
 Aalborg University, Fredrik Bajers Vej 7E, DK-9220 Aalborg, email:
\texttt{jm@math.aau.dk, rw@math.aau.dk}}

\maketitle

\begin{abstract}
This paper establishes a remarkable result regarding Palm
distributions for a log Gaussian Cox process: the reduced Palm
distribution for a log Gaussian Cox process is itself a log Gaussian Cox
process which only differs from the original log Gaussian Cox process
in the intensity function. This new result is used to study functional
summaries for log Gaussian Cox processes.
\end{abstract}


\section{Introduction}
Palm distributions  \citep[see e.g.][]{moeller:waagepetersen:04,daley:vere-jones:08}
are important in the theory and application of spatial
point processes. Intuitively speaking, for a prespecified location in space, the
Palm distribution of a point process, with respect to this location, plays the role of the conditional
distribution of the point process given that the aforementioned
location is occupied by a point of the point process.

The present paper focuses on  log Gaussian Cox processes 
\citep{moeller:syversveen:waagepetersen:98} which provide very
flexible, useful, and popular models for modeling spatial patterns in
e.g.\ biology
and spatial epidemiology. The paper establishes a
surprisingly simple characterization of
Palm distributions for such a process: The reduced $n$-point Palm distribution
is for any $n \ge 1$ itself  a log Gaussian Cox process that only differs from the
original log Gaussian Cox process in its intensity function (not to be confused with the random intensity function generating this kind of Cox process). This
result can be exploited for functional summaries as discussed
later. The simplicity and completeness of this result is
  remarkable when compared with Palm distributions for other common
  classes of spatial point processes. Reduced Palm distributions for Gibbs
  point processes are also themselves Gibbs point processes but with
  densities only known up to a normalizing constant. For shot-noise
  Cox processes \citep{moeller:02} one-point reduced Palm
  distributions have a simple characterization as cluster processes
  similar to shot-noise Cox processes
  but this is not the case for $n$-point Palm distributions when $n
  >1$.

The paper is organized as follows. 
Section~\ref{sec:palmandpre} reviews the general definition of
reduced Palm
distributions of any order and relates this to Cox processes. 
Section~\ref{s:lgcp} establishes our
characterization result for log Gaussian Cox processes.
Section~\ref{s:FGJ} applies this result to
functional summaries
for stationary log
Gaussian Cox processes, in particular the so-called $F$,
$G$, and $J$-functions, where we establish 
some new theoretical results, consider how to
calculate $F$,
$G$, and $J$ using Laplace approximations, and discuss
an application.
Section~\ref{s:concluding} concludes the paper.

\section{Palm distributions}\label{sec:palmandpre}

Our general setting is as follows. For ease of exposition
we view a point process as a
random locally finite subset $\mathbf X$ of a 
Borel set $S \subseteq\mathbb R^d$, $d\geq 1$; 
for measure theoretical details, see e.g.\ 
\citet{moeller:waagepetersen:04} or
\citet{daley:vere-jones:03}.
Denoting $\mathbf X_B=\mathbf X\cap B$
the restriction of $\mathbf X$ to a set $B\subseteq S$, 
local finiteness of $\mathbf X$ means that 
$\mathbf X_B$
is finite almost surely (a.s.) whenever
$B$ is bounded.
We denote ${\cal N}$  
the state space consisting of 
the locally
finite subsets (or point configurations) of $S$. 
We use the generic notation $h$ for an arbitrary
 non-negative measurable function defined on ${\cal N}$, $S^n$, or $S^n \times
 {\cal N}$ for $n= 1,2,\ldots$.
Furthermore,
$\mathcal B_0$ is the family of all bounded
Borel subsets of $S$. 
Finally, 
recall that the void probabilities 
$\mathrm P(\bX_K=\emptyset)$,
$K\subseteq S$ compact, 
 uniquely determine the distribution
of $\bX$.

\subsection{Factorial moment measures and Palm distributions}

This section provides the general definition of reduced Palm
distributions of any order. For finite point processes
specified by a density, a simpler and more explicit definition 
is available as reviewed 
in \cite{moeller:coeurjolly:waagepetersen:tutorial}.

For $n=1,2,\ldots$ and  $B_i \in {\mathcal B}_0$, 
the $n$-th order factorial moment measure
$\af^{(n)}$ is defined by
\[ \af^{(n)}(B_1 \times B_2 \times \cdots \times B_n) = \E \sum_{\xnl \in
    \bX}^{\neq} 1( x_1 \in B_1,\ldots,x_n \in B_n), \]
 where $1(\cdot )$ denotes the indicator function and $\neq$
over the summation sign means that $\xnl$ are pairwise distinct. 
If $\af^{(n)}$ has a density $\rho^{(n)}$ with respect to Lebesgue
measure, $\rho^{(n)}$ is called the $n$-th order joint intensity
function and is determined up to a Lebesgue nullset. Therefore, we can
assume that $\rho^{(n)}(x_1,\ldots,x_n)$ is invariant under
permutations of $\xnl$, and
 we need only to consider the case
where $x_1,\ldots,x_n\in S$ are pairwise distinct. Then 
$\rho^{(n)}(x_1,\ldots,x_n)\,\mathrm dx_1\cdots\,\mathrm dx_n$ can be
interpreted as the approximate probability for $\bX$ having a point in
each of infinitesimally small regions around $x_1,\ldots,x_n$ of
volumes $\mathrm dx_1,\ldots\,\mathrm dx_n$, respectively. We also write $\rho(u)$ for the intensity function $\rho^{(1)}(u)$.

Moreover, for any measurable $F \subseteq {\mathcal N}$,
define the $n$-th order reduced Campbell measure $C^{(n)!}$ as the
measure on $S^n\times\mathcal N$ given by
\begin{align*}
C^{(n)!}&( B_1 \times B_2 \times \cdots \times B_n \times F) =\\& \E \sum_{\xnl \in
    \bX}^{\neq} 1( x_1 \in B_1,\ldots,x_n \in B_n, \bX \setminus
  \{\xnl\} \in F ). 
\end{align*}
 Note that
$C^{(n)!}(\cdot \times F)$, as a measure on $S^n$,
is absolutely
continuous with respect to $\af^{(n)}$, with a density
$P^{!}_{\xnl}(F)$ which is determined up to an $\af^{(n)}$ nullset, and
$\af^{(n)}(\mathcal N)=C^{(n)!}( B_1 \times B_2 \times \cdots \times B_n \times \mathcal N)$. By the so-called Campbell-Mecke formula/theorem,
we can assume that $P^{!}_{\xnl}(\cdot)$ is a point process 
distribution on ${\mathcal N}$, called the $n$-th order reduced
Palm distribution given $\xnl$ \citep[see e.g.][]{daley:vere-jones:08}. 
We denote by $\bXred$ a point process
distributed according to $P^{!}_{\xnl}$. Again we need only to
consider the case where $\xnl$ are pairwise distinct. Then $P^{!}_{\xnl}$
can be interpreted as the conditional distribution of
$\bX\setminus\{\xnl\}$ given that $\xnl\in\bX$. 

If $\rho^{(n)}$ exists, then 
by standard measure theoretical arguments we obtain the extended Campbell-Mecke formula
\begin{align}
  &  \E \sum_{\xnl \in \bX}^{\neq} h(\xnl,\bX
\setminus \xn ) \nonumber \\
= &\, \int_{S}\cdots\int_{S} \E h(\xnl,\bXred ) \rho^{(n)}(\xnl) \dd x_1
\cdots \dd x_n . \label{eq:palmdef}
\end{align}
Suppose $\rho^{(m+n)}$ exists for an $m \ge
1$ and $n \ge 1$. Then,
for  pairwise distinct $u_1,\ldots,u_m,\xnl \in S$, 
it follows easily by expressing $\af^{(m+n)}$ as an
expectation of the form \eqref{eq:palmdef} that
$\mathbf X_{\xnl}^!$ has $m$-th order joint intensity function
\begin{equation}\label{eq:palmprodintens}
\rho^{(m)}_{x_1,\ldots,x_n}(u_1,\ldots,u_{m})=\left\{
\begin{array}{ll}
\frac{\rho^{(m+n)}(u_1,\ldots,u_{m},x_1,\ldots,x_n)}{\rho^{(n)}(x_1,\ldots,x_n)} & 
\mbox{if $\rho^{(n)}(x_1,\ldots,x_n)>0$,}\\
0 & \mbox{otherwise.}
\end{array}
\right.
\end{equation}
We also write $\rho_{x_1\dots,x_n}$ for the intensity function $\rho^{(1)}_{x_1\dots,x_n}$.

\subsection{Cox processes}

Let $\mathbf \Lambda=\{\Lambda(x)\}_{x\in S}$ be a 
nonnegative random field such that  $\mathbf \Lambda$ is locally
integrable a.s., that is,
for any $B\in\mathcal B_0$, the integral
$\int_B\Lambda(x)\,\mathrm dx$ exists and is finite a.s. In the
sequel, $\mathbf X$ conditional on $\mathbf
\Lambda$ is assumed to be a Poisson process with intensity function
$\mathbf \Lambda$; we say that $\mathbf X$ is
 a Cox process driven by~$\mathbf
\Lambda$. We also assume that $\mathbf \Lambda$ has moments
of any order $n=1,2,\ldots$. Then the joint intensities of $\mathbf X$
exist: For any $n=1,2,\ldots$ and pairwise distinct $x_1,\ldots,x_n\in S$, 
\begin{equation}\label{eq:coxprodintens}
\rho^{(n)}(x_1,\ldots,x_n)=\mathrm
E\left\{\prod_{i=1}^n\Lambda(x_i)\right\}.
\end{equation}

The following lemma, which is verified in Appendix~A, 
gives a characterization of the reduced Palm
distributions and their void probabilities.
\begin{lemma}\label{lm:voidcoxpalm}
Let $\bX$ be a Cox process satisfying the conditions above. Then, 
for any $n=1,2,\ldots$, pairwise distinct $x_1,\ldots,x_n\in S$, 
and compact $K \subseteq S$,
\begin{equation}\label{e:cox1}
\rho^{(n)}(x_1,\ldots,x_n)\E\left\{h\left(\xnl,\bXred\right)\right\}=
\E\left\{h(\xnl,\bX)\prod_{i=1}^n\Lambda(x_i)\right\}
\end{equation}
and
\begin{equation}\label{e:void12}
\rho^{(n)}(x_1,\ldots,x_n) 
\PP(\bXred \cap K = \emptyset ) = \E \left [ 
\exp\left\{ -\int_K \Lambda(u) \dd u \right\}
\prod_{i=1}^n \Lambda(x_i)  \right ]. \end{equation}
\end{lemma}

\section{Reduced Palm distributions for log Gaussian Cox
  processes} \label{s:lgcp}

For the remainder of this paper, let $\bX$ be a 
Cox process driven by $\mathbf \Lambda = \{\Lambda(x)\}_{x \in S}$,
 where $\Lambda(x)=\exp\{Y(x)\}$ and
$\mathbf Y=\{Y(x)\}_{x\in S}$ is a Gaussian random field with
mean function $\mu$ and covariance function $c$ 
so that 
$\mathbf\Lambda$ is locally integrable a.s.\ \citep[simple conditions ensuring this 
are given in][]{moeller:syversveen:waagepetersen:98}.
 Then $\mathbf X$ is a log Gaussian Cox
process (LGCP) as introduced by \citet{coles:jones:91} in astronomy  and
independently by \citet{moeller:syversveen:waagepetersen:98} in statistics.

For distinct $x,y\in S$, define the so-called pair correlation function $g(x,y)=\rho^{(2)}(x,y)/ \{ \rho(x)\rho(y)\}$ (the following result shows that $\rho>0$ in the present case).
By \citet[][Theorem~1]{moeller:syversveen:waagepetersen:98},
\begin{equation}\label{eq:pcflgcp}
\rho(x)=\exp\{\mu(x)+c(x,x)/2\},\qquad 
g(x,y)=\exp\{c(x,y)\},\qquad x,y\in S,\ x\not=y,
\end{equation}
and for pairwise distinct $x_1,\ldots,x_n\in S$, 
\begin{equation}\label{eq:prodintenslgcp}
\rho^{(n)}(x_1,\ldots,x_n)=\left\{\prod_{i=1}^n\rho(x_i)\right\}
\left\{\prod_{1\le i<j\le n}g(x_i,x_j)\right\}
\end{equation}
is strictly positive. 

For $u,x_1,\ldots,x_n\in S$, define
\[\mu_{\xnl}(u)=\mu(u)+\sum_{i=1}^nc(u,x_i).\] 
Combining \eqref{eq:palmprodintens} and
\eqref{eq:pcflgcp}-\eqref{eq:prodintenslgcp}, we obtain
for any pairwise distinct \linebreak
$u_1,\ldots,u_m,x_1,\ldots,x_{n}\in S$ with $m>0$ and $n>0$,
\begin{equation}
\rho^{(m)}_{x_1,\ldots,x_n}(u_{1},\ldots,u_{m})=
\left\{\prod_{i=1}^{m}\rho_{\xnl}(u_i)\right\}\left\{\prod_{1\le
    i<j\le m}g(u_i,u_j)\right\}
\label{e:9a},
\end{equation}
where 
\[\rho_{x_1,\ldots,x_n}(u)=
\exp\left\{\mu_{\xnl}(u)+c(u,u)/2\right\}.\]
Thereby the following proposition follows.

\begin{proposition}\label{p:2}
For the LGCP $\bX$ and any pairwise distinct $x_1,\ldots,x_n\in S$, 
$\bXred$ has $m$-th order joint intensity \eqref{e:9a} 
which agrees with the $m$-th order joint intensity
 function for an LGCP with mean function $\mu_{x_1,\ldots,x_n}$
 and covariance function $c$ for the underlying Gaussian random field.
\end{proposition}

Proposition~\ref{p:2} indicates that also $\bXred$ could be an LGCP.
A sufficient condition, considered by \cite{macchi:75}, is the existence of a number $a=a(B)>1$ for
each set $B\in \mathcal B_0$  such that
\begin{equation}\label{eq:conditionMacchi}
  \mathrm E \left[ \exp \left\{ a \int_B \Lambda(u) \dd u\right\}\right] <\infty.
\end{equation}
 However,
 we have not been
successful in verifying this condition which
 seems too strong to hold for
any of the covariance function models we have considered, including when
$c$ is constant (then $\bX$
is a mixed Poisson process) or weaker cases of correlation,  
e.g.\ if $c$ is a stationary
exponential covariance function. The case where
$c$ is constant is closely related to the log normal
  distribution which is not
uniquely determined by its moments \citep{heyde:63}. 

Accordingly we use instead 
Lemma~\ref{lm:voidcoxpalm} when establishing 
the following theorem, which implies that the LGCPs $\bX$ and
$\bXred$ share the same pair correlation function and differ only in
their intensity functions.
\begin{theorem}\label{thm} For pairwise distinct $\xnl\in S$, 
$\bXred$ is an LGCP with  underlying
  Gaussian random field $\bY_\xnl$, where $\bY_\xnl$ has mean function $\mu_{\xnl}$ and 
covariance function $c$. 
\end{theorem}

Let ${\tilde{\bY}}=\bY-\mu$ be the centered Gaussian random field with
covariance function~$c$. 
Theorem~\ref{thm} is a consequence of the fact that the probability measure of $\bY_\xnl$ is
absolutely continuous with respect to the one of ${\tilde{\bY}}$,
 with density $\exp \left\{
  \sum_{i=1}^n {\tilde{y}}(x_i)- \sum_{i,j=1}^n c(x_i,x_j)/2 \right\}$
when ${\tilde{\by}}$ is a realization of ${\tilde{\bY}}$. This
result is related to the Cameron-Martin-Girsanov formula for
one-dimensional Gaussian processes. A short selfcontained proof
covering our spatial setting is given in Appendix~A.

Often we consider a non-negative covariance function $c$ or equivalently
 $g\ge1$, which is interpreted as `attractiveness of the LGCP at
all ranges', but even more can be said: 
A coupling between $\bX$ and $\bXred$ is obtained by taking
$Y_{x_1,\ldots,x_n}(x)=Y(x)+\sum_{i=1}^n c(x,x_i)$. 
Thus, if $c\ge0$ and we are given pairwise distinct points $x_1,\ldots,x_n \in S$, we can consider ${\mathbf
X}$ as being included in ${\mathbf
X}_{x_1,\ldots,x_n}^!$, since ${\mathbf
X}$ can be obtained by an independent thinning
of $\mathbf
X_{x_1,\ldots,x_n}^!$, with inclusion probabilities 
$\exp\{-\sum_{i=1}^n c(x,x_i)\}$, $x\in \bX\setminus\{x_1,\ldots,x_n\}$.
This property clearly shows the attractiveness of the LGCP if $c\ge0$
(equivalently $g\ge1$).

\section{Functional summaries for stationary  log Gaussian Cox
  processes }\label{s:FGJ}

Throughout this section, 
let $S=\mathbb R^d$ 
and assume that the LGCP $\bX$ is stationary, i.e., its distribution
is invariant under translations in $\mathbb R^d$.
By 
\eqref{eq:pcflgcp}-\eqref{eq:prodintenslgcp}, 
this is equivalent to 
stationarity of the 
underlying Gaussian random field $\bY$, that is, 
the intensity $\rho$ is
constant and the  pair correlation function $g(x,y)=\giso (x-y)$ is translation
invariant, where $x,y\in \R^d$ are distinct, and
$\giso(x)=\exp\{ \ciso(x)\}$ and $\ciso(x)=c(o,x)$ for $x\in\mathbb R^d$, where $o$ denotes the
origin in $\mathbb R^d$. It is custom to call $\mathrm P^!_o$ the
reduced Palm distribution at a typical point, noticing that for any
$x\in \mathbb R^d$, $\bX^!_o$ and $\bX^!_x-x=\{y-x:y\in\bX^!_x\}$ 
are identically distributed. 

Denote $B(o,r)$ the
ball in $\mathbb R^d$ of radius $r>0$ and centered at $o$. 
Popular tools for exploratory purposes as well as
model fitting and model checking are based on the following functional 
summaries where $r>0$ \citep[see e.g.][]{moeller:waagepetersen:04}:
\begin{itemize}
   \item[(i)] 
   the pair correlation function
$\giso$ 
and
the related Ripley's $K$-function given by 
\[K(r)=\frac{1}{\rho}\mathrm E\,\#\left\{\bX^!_o\cap B(o,r)\right\}=
\int_{B(o,r)}\giso (x).\]
Thus, $\rho K(r)$ is the expected number of further points in $\bX$ within
distance $r$ of a typical
point in $\bX$. If $\giso(x)$ depends only on the distance $\|x\|$
then $\giso$ and $K$ are in one-to-one correspondence;
\item[(ii)] the empty space function given by
\[F(r) = \mathrm P \left\{ \bX\cap B(o,r)\not=\emptyset\right\},\]
which is the probability that $\bX$ has a point within distance $r$
of an arbitrary fixed location;
\item[(iii)] the nearest-neighbour distribution function given by
\[G(r) = \mathrm P \left\{\bX^!_o\cap B(o,r)\not=\emptyset \right\},\]
which is the probability that $\bX$ has a further point within
distance $r$ of a typical
point in $\bX$;
\item[(iv)] the $J$-function given by
\[J(r) = \frac{1-G(r)}{1-F(r)},\] 
with the convention $a/0=0$ for any $a\geq0$.
\end{itemize} 

Section~\ref{s:new}
establishes some new results for these theoretical functions and 
Section~\ref{s:laplace} discusses how they can be calculated using a
Laplace approximation. Section~\ref{s:numerical} illustrates this
calculation and Section~\ref{s:application} discusses an application
for a real dataset. 

\subsection{New formulae for $G$~and~$J$}\label{s:new}

By conditioning on $\bY$, we see that 
\begin{equation}
1-F(r) = \mathrm E \left(
\exp \left[ -\int_{B(o,r)} \exp\{ Y(x)\} \dd x
\right]
\right).\label{eq:F1} 
\end{equation}
Using the 
Slivnyak-Mecke formula, 
\citet{moeller:syversveen:waagepetersen:98} showed that
\begin{equation}
1-G(r) 
= \frac1\rho \; \mathrm E \left(
\exp \left[ Y(o)-\int_{B(o,r)} \exp\{ Y(x)\} \dd x
\right]
\right) . \label{eq:G1}
\end{equation}
Since the nearest-neighbour distribution function for $\bX$ is the same as the 
empty space function for $\bX^!_o$, which is an LGCP
with underlying Gaussian random field
$Y_o(x)=Y(x)+\ciso(x)$, and since $\giso(x)=\exp\{\ciso(x)\}$, 
we obtain an alternative expression
\begin{equation}
1-G(r)=\mathrm E \left(
\exp \left[ -\int_{B(o,r)} \giso(x)\exp\{ Y(x)\} \dd x
\right]
\right).  \label{eq:G2}
\end{equation}
Therefore, we also obtain a  new expression for the $J$-function,
\begin{equation} \label{eq:newJ}
  J(r)  = \frac{\mathrm E \left(
\exp \left[ -\int_{B(o,r)} \giso(x)\exp\{ Y(x)\} \dd x
\right]
\right)}{
  \mathrm E \left(
\exp \left[ -\int_{B(o,r)} \exp\{ Y(x)\} \dd x
\right]
\right)
}.
\end{equation}

\citet{vanlieshout:11} established for a general stationary point process
the approximation $J(r)-1 \approx -\rho\{ K(r)-\omega_d r^d\}$,
where $\omega_d=|B(o,1)|$ and $r\mapsto \omega_d r^d$ is the $K$-function for a stationary Poisson
process. It is therefore not so surprising that often empirical 
$J$ and $K$-functions lead to the same practical interpretations. 
In particular, if for our LGCP $\ciso\ge0$, i.e., $\giso\geq 1$, 
then we have $K(r)-\omega_d r^d\geq 0$, and so we expect that $J(r)\leq
1$. Indeed \citet{vanlieshout:11} verified this in the case of an
LGCP with $\giso\geq 1$. This result immediately follows by
the new expression~\eqref{eq:newJ}.

\subsection{Laplace approximation}\label{s:laplace}

Since Laplace's pioneering work \citep[see e.g.][]{stigler:86}, Laplace approximations of complex integrals have gained much attention in probability and statistics, 
in particular when considering integrals involving Gaussian random fields \citep[see e.g.][]{rue:martino:chopin:09}.  This section discusses a Laplace approximation of $1-G(r)$; a Laplace
approximation of $1-F(r)$ can be obtained along similar lines. 

For $\Delta>0$, 
consider a grid of quadrature points,
\[ {\cal G}(\Delta)=\{ (\Delta
i_1,\dots, \Delta i_d) \, | \, i_1,\dots,i_d \in \mathbb{Z} \},\] 
and for $v \in
{\cal G}(\Delta)$, 
let $A^\Delta_v=[v_1-\Delta/2,v_1+\Delta/2 [ \times \dots \times
[v_d-\Delta/2,v_d+\Delta/2[$ be the grid cell associated with $v$.
Then for any non-negative Borel function $\ell:\R^d\to \R$, 
we use the numerical quadrature approximation
\begin{equation}\label{e:app0}
  \int_{B(o,r)} \exp\{ Y(x)\} \ell(x) \dd x \approx \sum_{v\in
          {\cal G}(\Delta) \cap B(o,r) } w_v \ell(v) \exp\{Y(v)\},
\end{equation}
where the quadrature weight $w_v= |A_v^\Delta \cap B(o,r)|$. 

Denote by $f$, $M$ and $\Sigma$ the density, the mean
vector and the covariance matrix of the normally distributed vector 
$\{Y(v)\}_{v\in {\cal G}(\Delta) \cap B(o,r)}$. 
Then \eqref{eq:G2} and \eqref{e:app0} give
\begin{align}
  1-G(r) & \approx \int_{\R^m} \!\!\! \exp \left \{- \!\!\!\!\!\! \sum_{v\in {\cal
        G}(\Delta) \cap B(o,r)}  \!\!\!\!\!\! w_v \giso(v)\exp(y_v) \right \}
  f(y) \dd y \nonumber \\ &  = \int_{\R^m} \!\!\! \exp \{h(y) \} \dd y \label{e:app1}
\end{align}
where $y$ is the vector $(y_v)_{v\in {\cal G}(\Delta) \cap B(o,r)}$ of 
dimension $m=\#\{\mathcal G(\Delta) \cap B(o,r)\}$ and 
\[h(y)= - \!\!\!\!\!\! \sum_{v\in {\cal G}(\Delta) \cap B(o,r)}  \!\!\!\!\!\! w_v \giso(v)\exp(y_v)-\frac12
(y-M)^\top\Sigma^{-1}(y-M)-\frac12 \log\{ (2\pi)^{m} |\Sigma|
\}.\] 
The gradient vector for $h$ is
\begin{equation} \label{e:sigma}
  \nabla h(y) = -d(y) - \Sigma^{-1}(y-M),
\end{equation}
where $d(y)=\{w_v \giso(v) \exp(y_v)\}_{v\in {\cal G}(\Delta) \cap
  B(o,r)}$, and minus one times the Hessian matrix for $h$ is
\[
  \quad H(y) = D(y)+\Sigma^{-1},
\]
where  $D(y)$ is the diagonal matrix with entries $\{d(y)\}_v$, $v  {\cal G}(\Delta) \cap B(o,r) \}$. 
Since $H(y)$ is a positive definite matrix, $h$ has a unique maximum at a point
$\hat y$, which can be found using Newton-Raphson iterations 
\begin{equation} \label{e:yH}
y^{(l+1)}=y^{(l)}+ H^{-1}\{y^{(l)} \}\nabla h\{y^{(l)}\}. 
\end{equation}
Therefore, 
the logarithm of the Laplace approximation of the right hand
side in \eqref{e:app1} \citep[see e.g.][]{stigler:86}
gives
\begin{align} 
  \log \{ 1-G(r) \} \approx& -\sum_{v\in {\cal G}(\Delta) \cap B(o,r)}  w_v \giso(v) \exp(\hat y_v)  +\frac12 (\hat y-M)^\top d(\hat y) \nonumber\\ 
  &-\frac12 \log| D(\hat y)\Sigma+I|.\label{eq:Glapl}
\end{align}
where $I$ is the 
$m\times m$ identity matrix.
For the computation of $\Sigma^{-1}(y-M)$ in \eqref{e:sigma}
we solve $LL^\top z=y-M$ where $L$ is the Cholesky factor of
$\Sigma$. 
In the same way, considering the $QR$ decomposition, $Q(y)R(y)$ for $y\in \R^m$, of the matrix 
$D(y)\Sigma+I$, the computation of $H^{-(1)}\{y^{(l)}\} \nabla
h\{y^{(l)}\}$ in \eqref{e:yH} is done by first solving
$Q\{y^{(l)}\}R\{{y^{(l)}}\} \tilde z=\nabla h\{y^{(l)}\}$ and second
by evaluating $\Sigma \tilde z$. Finally, 
in \eqref{eq:Glapl}, $|D(\hat y)\Sigma+I|= |R(\hat y)|$.

\subsection{Numerical illustration}\label{s:numerical}

\begin{figure}[htbp]
\centering
\hspace*{-.5cm}\begin{tabular}{ll}
  \includegraphics[scale=.45]{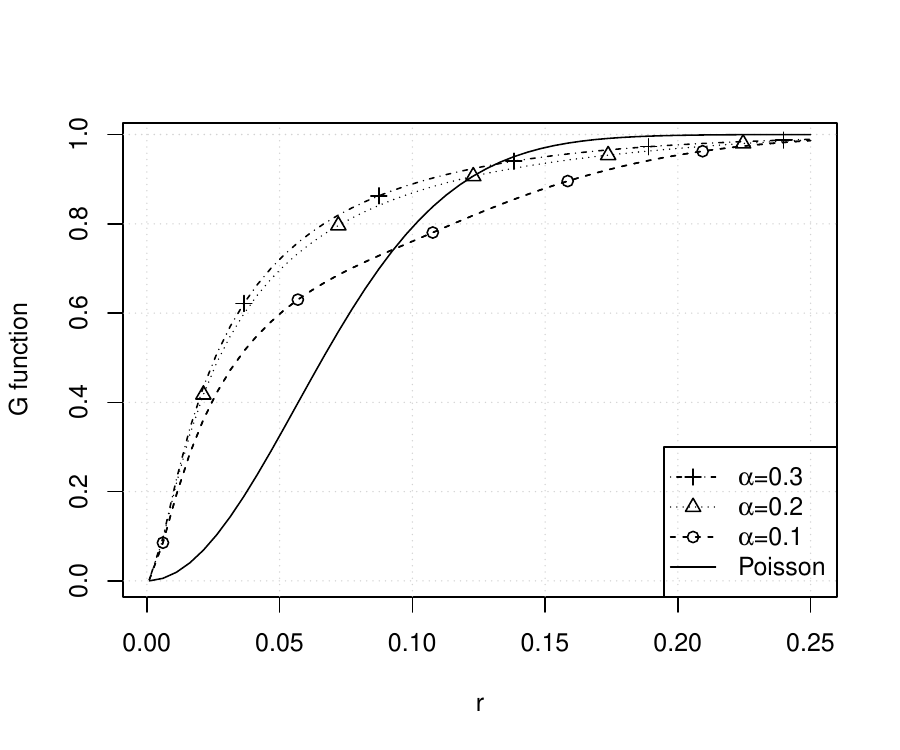} &
  \includegraphics[scale=.45]{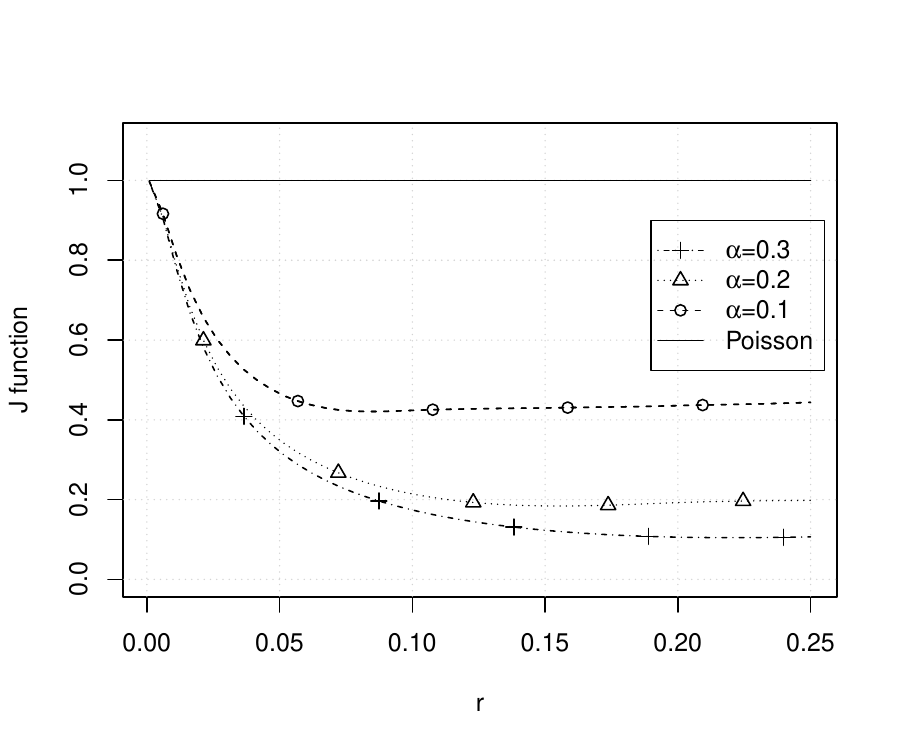} \vspace*{-.2cm}
\end{tabular}
  \caption{\label{fig:lapl} 
For three planar stationary LGCPs with intensity $\rho=50$ and $\tilde
c$ given by a
       spherical covariance function, with variance $\sigma^2=4$ and
       scale parameters $\alpha=0.1,0.2,0.3$, respectively,
Laplace approximations of the 
$G$-function (left) and the $J$-function (right). 
}
\end{figure}

To illustrate the Laplace approximations of the $G$ and $J$-functions (Section~\ref{s:laplace})
we consider three planar stationary LGCPs
 with intensity $\rho=50$ and spherical covariance function 
\[ \tilde c(x)=
\left\{ 
\begin{array}{ll}
\sigma^2\left[1-\frac2\pi\left\{\frac{\|x\|}{\alpha}
{\sqrt{1-\left(\frac{\|x\|}{\alpha}\right)^2}}+\sin^{-1}
\frac{\|x\|}{\alpha}\right\}\right] & 
\mbox{if } \|x\|\leq \alpha, \\
0 & \mbox{otherwise,} 
\end{array}
\right. 
\]
with
 variance $\sigma^2=4$ and scale parameters $\alpha=0.1,0.2,0.3$,
 respectively. We evaluate  the approximations of $G(r)$ and
 $J(r)$ 
 at $r\in \mathcal R$, where $\mathcal R$ is the set of 50 equispaced values between 0.01 and
 0.25. 
For $r \in \mathcal R$, we define the grid $\mathcal G(\Delta_r)$
  with $\Delta_r=2r/q$, where $q$ is a fixed integer. Such a choice
  implies that $\#\{\mathcal G(\Delta_r)\cap [-r,r]^2\}=q^2$, and so we
  have at least $q^2 \pi/4$ quadrature points in $B(o,r)$.  For a
given $q$, denote by $G_q$, $F_q$, and $J_q$ the corresponding
Laplace approximations of $G$, $F$, and $J$, respectively.
Figure~\ref{fig:lapl} shows the resulting curves with
 $q=16$. 
To see how far these Cox processes deviate from the Poisson case (which would correspond to $\sigma^2=0$), we also plot the $G$-function in the Poisson case, namely $1-G(r)=\exp(-\rho \pi r^2)$.
To study the role of $q$, we 
report in Table~\ref{tab1} the maximal differences $\max_{r\in
  \mathcal R}|G_{16}(r)-G_q(r)|$ and $\max_{r\in \mathcal R} |J_{16}(r)-J_q(r)|$ for $q=4,8,12$.
As expected, each difference decreases as $q$ increases and is already very
small when $q=12$ (less than $4\times 10^{-3}$ except for the
$J$-function and $\alpha=1$). This justifies our choice $q=16$ in
Figure~\ref{fig:lapl}. 

 \begin{table}[H]
\centering
\begin{tabular}{rrrr}
  \hline
 & \qquad $q=4$ &\qquad  $q=8$ &\qquad $q=12$ \\ 
  \hline
$\alpha=0.1$, $G$ & 59.9 & 8.4 & 2.1 \\
  $J$ & 505.9 & 96.1 & 20.5 \\
$\alpha=0.2$, $G$  & 14.3 & 1.6 & 0.5 \\
  $J$ &  109.0 & 13.8 & 3.5 \\
  $\alpha=0.3$, $G$ &4.2 & 0.5 & 0.1 \\
  $J$ & 22.1 & 3.1 & 0.3 \\
   \hline
\end{tabular}
\caption{\label{tab1} For the same three LGCPs
  as in Figure~\ref{fig:lapl}, maximal differences between the Laplace
  approximations $G_q$ and $G_{16}$, and between 
the Laplace
  approximations $J_q$ and $J_{16}$, with
  $q=4,8,12$, that is, the evaluation of $\max_{r\in \mathcal R}|H_{16}(r)-H_q(r)|$ for $H=G,J$. Results are multiplied by $10^3$.}
\end{table}

The Laplace approximation of the $G$-function could also be derived
using \eqref{eq:G1}. To check the
  agreement of the numerical approximations based on
\eqref{eq:G1} and~\eqref{eq:G2}, respectively,
Table~\ref{tab2} shows the maximal difference between the two approximations
of first the $G$-function and second the $J$-function. In agreement with the
theoretical developments, in both cases,
the difference does not exceed $4\times 10^{-4}$ when $q=16$. 

\begin{table}[ht]
\centering
\begin{tabular}{rrrrr}
  \hline
 & \quad$q=4$ &\quad  $q=8$ &\quad $q=12$ &\quad $q=16$ \\
  \hline
$\alpha=0.1$, $G$ &  3.8 & 8.4 & 4.4 & 3.2 \\
  $J$ &15 & 8.5 & 6.1 & 3.9 \\
  $\alpha=0.2$, $G$ & 4.7 & 3.1 & 3.9 & 1.5 \\
  $J$ &  4.7 & 3.1 & 2.1 & 1.7 \\
  $\alpha=0.3$, $G$ &  0.1 & 1.9 & 1.4 & 1.1 \\
  $J$ &  0.2 & 2.0 & 1.5 & 1.3 \\ 
   \hline
\end{tabular}
\caption{\label{tab2}  For the same three LGCPs
  as in Figure~\ref{fig:lapl},
maximal differences between 
the Laplace
  approximations of the $G$-function based on 
\eqref{eq:G1} respective
  \eqref{eq:G2}, with $q=4,8,12,16$, and similarly 
for the $J$-function. Results are multiplied by $10^4$.}
\end{table}

\subsection{Scots pine saplings dataset}\label{s:application}

\begin{figure}[H]
\hspace*{-2cm}  \begin{tabular}{ll}
    \includegraphics[scale=.47]{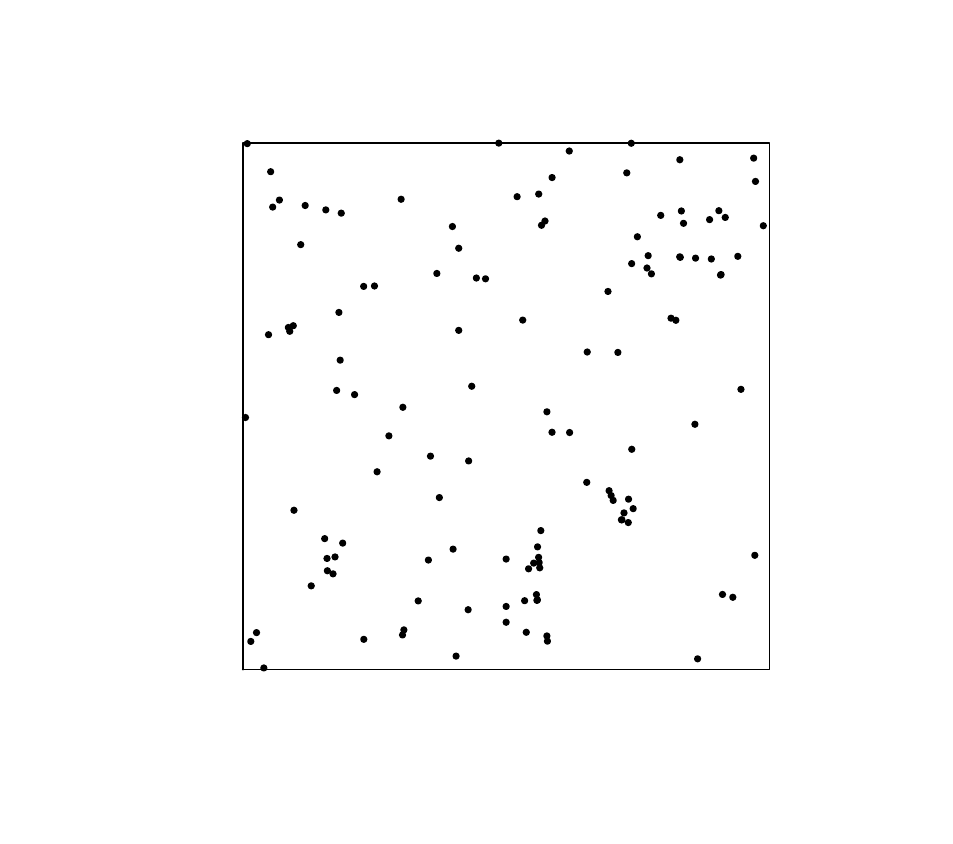} &
    \includegraphics[scale=.51]{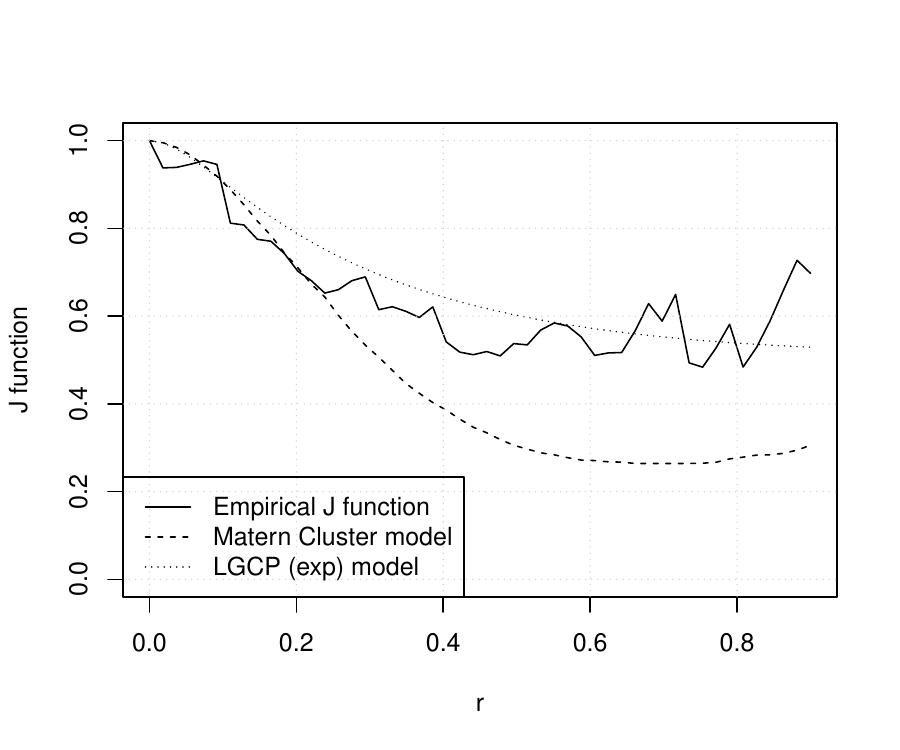} \vspace*{-.5cm}
  \end{tabular}
  \caption{\label{fig:pines} Left panel: Locations of 126 Scots
          pine saplings in a 10 by 10 metre square. Right panel:
          Non-parametric estimate of the $J$-function (solid curve) and
          fitted $J$-functions for the Mat{\'e}rn cluster process
          (dashed curve) and the LGCP with exponential covariance
          function (dotted curve).}
\end{figure}

The left panel in Figure~\ref{fig:pines} shows 
the locations
of 126 Scots pine saplings in a 10 by 10 metre square. 
The dataset is included in the 
\texttt{R} package \texttt{spatstat} as \texttt{finpines}, and
it has previously been
analyzed by \citet{penttinen:stoyan:henttonen:92},
\citet{stoyan:stoyan:94}, and
\citet{moeller:syversveen:waagepetersen:98}. In the first two 
papers a Mat{\'e}rn cluster process is fitted, using the $K$-function (or
its equivalent $L$-function) and its nonparametric estimate 
both for parameter estimation and model
checking, while the third paper considered an LGCP with exponential
covariance function and used the pair correlation function for
parameter estimation and the $F$ and $G$-functions for model
checking. \citet{moeller:syversveen:waagepetersen:98} concluded that
both models provide a reasonable fit although when also including
a third-order functional summary (i.e., one based on
$\bX^!_{o,x}$)  
the LGCP model showed a better fit. Below we extend this analysis by
using the $J$-function and 
the approximation established in Section~\ref{s:laplace}.

We fitted both models by minimum contrast estimation (method
\texttt{kppm} in
\texttt{spatstat}) which compares a non-parametric estimate of the 
$K$-function with its theoretical value. 
When approximating the $J$-function for the LGCP, we used 
$q=12$; no improvements were noticed with higher values of
  $q$. 
The right panel in Figure~\ref{fig:pines} shows the theoretical
$J$-functions for
the two fitted models together with a non-parametric estimate of the 
$J$-function obtained from data, considering 50 equispaced distances ($r$-values)  between 0 and 0.9 meter
\citep[for the exact expression of
the $J$-function for the Mat{\'e}rn cluster process, see e.g.][]
{moeller:waagepetersen:04}.
Clearly,  
the fitted LGCP provides a better fit than the fitted
Mat{\'e}rn cluster process. Indeed, 
the maximal difference between the non-parametric estimate
and the theoretical $J$-function equals $0.43$ for the Mat{\'e}rn cluster model and $0.20$ for the LGCP model.

\section{Concluding remarks}\label{s:concluding}

We expect that our results for 
the reduced Palm distributions for an LGCP can be exploited further
 regarding
third-order and higher order functional
summaries (one such characteristic was briefly
mentioned in Section~\ref{s:application}), parameter estimation
procedures, model checking, etc. This is discussed in some detail below.

For likelihood based inference, suppressing in the notation any
dependence on a parametric model and assuming a realization $\mathbf
x$ of an LGCP is observed within a region $W\subset\mathbb R^d$ of Lebesgue measure $|W|\in(0,\infty)$, the likelihood function is given by the density
\[
  f(\mathbf x) = \E \left(
  \exp\left[ 
|W|-\int_W \exp \{ Y(u)\} \dd u
  \right] 
  \prod_{u\in \mathbf x} \exp\{ Y(u)\}
\right)
\]
with respect to the unit rate Poisson process on $W$. 
This density expression has no explicit form even for simple covariance function models for the underlying Gaussian random field. For maximum likelihood estimation and prediction of the Gaussian random field, rather elaborate and time-consuming Markov chain Monte Carlo (MCMC) methods \citep{moeller:waagepetersen:04} may be used; alternatively, a Bayesian approach based on integrated nested Laplace approximations or MCMC methods may be used \citep{rue:martino:chopin:09,taylor:diggle:14}. Using the results of this paper, we have the following new expression for the density:
\[
  f(\mathbf x) = \rho^{(n)} (x_1,\dots,x_n)  \E 
\left(
  \exp \left[
|W|-\int_W \exp \{ Y_{x_1,\dots,x_n}(u)\} \dd u 
  \right]\right). 
\]
Here the expression for the $n$-th order intensity $\rho^{(n)}$ is explicit, but it remains to investigate if the expectation, which apparently has a simpler expression, may be easier to approximate.

For model checking, when considering non-parametric estimates of functional summaries together with simulated confidence bounds under a fitted LGCP model \citep[such as extreme rank envelopes, see][]{baddeley:rubak:turner:15,myllymaki:etal:16}, 
  it could be pertinent to include the theoretical expressions of the functional summaries for LGCPs obtained in the present paper.

Finally, recall that for any point process,
the pair correlation
function (when it exists) is invariant under independent thinning. 
Could this property be 
exploited in connection to LGCPs where we know how the 
pair correlation
function is related to those of the reduced Palm distributions?  


\subsection*{Acknowledgments}
We thank the two reviewers, the associate editor and the editor for their careful
reading of our paper and their useful comments and suggestions which
helped us to improve the paper. 

J. M\o{}ller and R. Waagepetersen were supported by the Danish Council for Independent Research | Natural
Sciences,
grant 12-124675,
"Mathematical and Statistical Analysis of Spatial Data", and by
the "Centre for Stochastic Geometry and Advanced Bioimaging",
funded by grant 8721 from the Villum Foundation. J.-F. Coeurjolly was supported by ANR-11-LABX-0025 PERSYVAL-Lab (2011, project OculoNimbus).

\appendix

\section{Proofs}

{\it Proof of Lemma~\ref{lm:voidcoxpalm}:} By conditioning on $\mathbf \Lambda$, 
\eqref{eq:palmdef}
becomes
\begin{align}
& \E \E \left\{ \sum^{\not=}_{x_1,\ldots,x_n \in \bX}
h(x_1,\ldots,x_n,\bX \setminus \{x_1,\ldots,x_n\})\,\bigg |\, \mathbf \Lambda 
\right\}\nonumber\\
=  &\, \E \E\left\{\int_S\cdots\int_S 
h(x_1,\ldots,x_n,\bX)\prod_{i=1}^n\Lambda(x_i) \dd x_1\cdots \dd x_n
\,\bigg |\,  \mathbf\Lambda 
\right\}\label{e:bla1}\\
=  &\,
\int_S\cdots\int_S \E\left\{
h(x_1,\ldots,x_n,\bX)\prod_{i=1}^n\Lambda(x_i)\right\}\dd x_1\cdots
\dd x_n.
\label{e:bla11}
\end{align}
Here, in \eqref{e:bla1} we use 
that $\bX$ given $\mathbf \Lambda$ is a Poisson
process and apply the extended Slivnyak-Mecke theorem
\citep{moeller:waagepetersen:04}, and in \eqref{e:bla11} we use Fubini's
theorem. Combining \eqref{eq:palmdef} and 
\eqref{e:bla11}, we deduce \eqref{e:cox1}. Finally, \eqref{e:void12}
follows from \eqref{e:cox1} with $h(\xnl,\bfx)=1(\bfx \cap K= \emptyset)$.

\vspace{0.5cm}

\noindent{\it Proof of Theorem~\ref{thm}:} 
By \eqref{e:void12} and \eqref{eq:pcflgcp}-\eqref{eq:prodintenslgcp}, 
we just have to show that for any compact $K\subseteq S$ and pairwise distinct
points $x_1,\ldots,x_n\in S$,
\begin{align*}
 &\E \exp\left[ -\int_K
   \exp\left\{{\tilde{Y}}(u)+\mu_\xnl(u)\right\}\dd u\right]\\
=\,&
 \E \exp \left [ \sum_{i=1}^n {\tilde{Y}}(x_i)- \sum_{i,j=1}^n c(x_i,x_j)/2 -\int_K
 \exp\left\{\mu(u)+{\tilde{Y}}(u)\right\} 
 \dd u \right ] .
\end{align*}
This will follow by verifying that the distribution of
$\{{\tilde{Y}}(u)+\sum_{i=1}^n c(u,x_i)\}_{u \in S}$ is absolutely continuous with
respect to the distribution of ${\tilde{\bY}}=\{{\tilde{Y}}(u)\}_{u \in S}$, with density $\exp \left\{
  \sum_{i=1}^n {\tilde{y}}(x_i)- \sum_{i,j=1}^n c(x_i,x_j)/2 \right\}$
when ${\tilde{\by}}$ is a realization of~${\tilde{\bY}}$. Since the
  distribution of a random field is determined by its finite
  dimensional distributions, we just need to verify 
the agreement of the
  characteristic functions of the {probability} measures $Q_1$ and $Q_2$
  given by 
\[ Q_1(B)=P \left[ \bigg\{ {\tilde{Y}}(u)+\sum_{i=1}^n c(u,x_i)
  \bigg\}_{u \in U} \in B  \right]\]
and
\[ Q_2(B)=\E \left( 1\left[\{ {\tilde{Y}}(u)\}_{u \in U}
\in B \right] \exp  \left\{\sum_{i=1}^n {\tilde{Y}}(x_i)- \sum_{i,j=1}^n c(x_i,x_j)/2 \right\}
\right), \]
for any Borel set $B\subseteq
\R^{n+m}$, any pairwise distinct
locations $u_1,\ldots,u_m \in \R^d \setminus \{ \xnl \}$, with
$m > 0$ and $U=\{x_1,\ldots,x_n,u_1,\ldots,u_m\}$. 
Let $\Sigma = \{c(u,v)\}_{u,v \in U}$ denote the covariance
matrix of  $\{{\tilde{Y}}(u)\}_{u \in U}$, and let
$c_\xnl=\{ \sum_{i=1}^n c(u,x_i)\}_{u \in U} = \Sigma e$, where
$e$ consists of $n$ 1's followed by $m$ 0's. For $t \in
\R^{n+m}$, the characteristic function of $Q_2$ is (with $i^2=-1$)
\begin{align*}  
& \E \exp \left[ i \{ {\tilde{Y}}(u)\}_{u \in I}^\T t     
+ \{ {\tilde{Y}}(u)\}_{u \in I}^\T e  + e^\T \Sigma e/2 \right] \\ 
&= \exp\left(e^\T \Sigma e/2\right) 
\E \exp \left[i \{ {\tilde{Y}}(u)\}_{u \in I}^\T (t - ie) \right] \\ 
&= 
\exp\left(e^\T \Sigma e/2 +i
e^\T \Sigma t - t^\T \Sigma t/2 - e^\T \Sigma e/2\right) \\ 
&=  \exp\left(i e^\T \Sigma t - t^\T \Sigma t/2\right). 
\end{align*}
The last expression is the characteristic function of $Q_1$ which
concludes the proof. For the second last equality in the above
derivation we considered $\{ {\tilde{Y}}(u)\}_{u \in U}$ as a complex
Gaussian vector and used the expression for its characteristic
function.


\bibliographystyle{plainnat}  
\bibliography{palm}

\end{document}